\documentclass{amsart}

\usepackage{amsmath,amscd,amssymb,amsthm}
\usepackage{latexsym,enumerate,graphicx,geometry,hyperref,rotating,braket}
\usepackage[all]{xy}

\newcommand{\thefilename}{ax2+by2=alpha}

\numberwithin{equation}{section}

\theoremstyle{plain}
\newtheorem{thm_}[equation]{Theorem}
\newtheorem{lemma_}[equation]{Lemma}
\newtheorem{prop_}[equation]{Proposition}
\newtheorem{cor_}[equation]{Corollary}
\newtheorem{eg_}[equation]{Example}

\theoremstyle{definition}
\newtheorem{thmu_}[equation]{Theorem}
\newtheorem*{thmus_}{Theorem}
\newtheorem{propu_}[equation]{Proposition}
\newtheorem{egu_}[equation]{Example}
\newtheorem{def_}[equation]{Definition}

\theoremstyle{remark}
\newtheorem{rk_}[equation]{Remark}

\newcommand{\thm}[1]{\begin{thm_}#1\end{thm_}}

\newcommand{\thmus}[1]{\begin{thmus_}#1\end{thmus_}}
\newcommand{\lemm}[1]{\begin{lemma_}#1\end{lemma_}}
\newcommand{\eg}[1]{\begin{eg_}#1\end{eg_}}

\newcommand{\rk}[1]{\begin{rk_}#1\end{rk_}}

\newcommand{\pf}[1]{\begin{proof}#1\end{proof}}
\DeclareMathOperator{\Gal}{Gal}

\newcommand{\fracl}[3]{\genfrac{(}{)}{}{}{#1}{#2}_{#3}}
\newcommand{\fracn}[2]{\genfrac{(}{)}{}{}{#1}{#2}}

\newcommand{\ZZ}{\mathbb Z}
\newcommand{\QQ}{\mathbb Q}
\newcommand{\RR}{\mathbb R}

\renewcommand{\AA}{\mathbb A}
\newcommand{\II}{\mathbb I}
\newcommand{\GG}{\mathbb G}
\newcommand{\XX}{\mathbf X}
\newcommand{\TT}{\mathbf T}

\renewcommand{\o}{\mathfrak o}
\newcommand{\p}{\mathfrak p}
\renewcommand{\P}{\mathfrak P}

\renewcommand{\t}{\times}

\newcommand{\lr}{\longrightarrow}

\newcommand{\eq}[1]{\begin{equation}#1\end{equation}}
\newcommand{\eqn}[1]{\begin{equation*}#1\end{equation*}}

\newcommand{\aln}[1]{\begin{align*}#1\end{align*}}

\newcommand{\cs}[1]{\begin{cases}#1\end{cases}}

\newcommand{\enmt}[1]{\begin{enumerate}#1\end{enumerate}}
\renewcommand{\it}{\item}

\newcommand{\itm}[1]{\it[\upshape{(#1)}]}
\newcommand{\newnoindbf}[1]{\vspace{2mm}\noindent\textbf#1}

\begin{document}
\title[Integral Representation of $ax^2+by^2$]{On the Integral Representation of $ax^2+by^2$ and the Artin Condition}
\author[C. Lv]{Chang Lv}
\address{Key Laboratory of Mathematics Mechanization\\
NCMIS, Academy of Mathematics and Systems Science\\
Chinese Academy of Sciences, Beijing 100190, P.R. China}
\email{lvchang@amss.ac.cn}
\author[J. Shentu]{Junchao Shentu}
\address{Institute of Mathematics\\
Academy of Mathematics and Systems Science\\
Chinese Academy of Sciences, Beijing 100190, P.R. China}
\email{stjc@amss.ac.cn}
\date{\today}
\begin{abstract}
Given a number field $F$ with $\o_F$ its ring of integers. For certain $a,b$ and  $\alpha$ in $\o_F$, we show that 
 the Artin condition is the only obstruction to the local-global principle for integral solutions of equation $ax^2+by^2=\alpha$. 
Some concrete examples are presented at last.
\end{abstract}
\maketitle

\section{Introduction}\label{sec_intro}
The main theorem of David A. Cox \cite{cox} is a beautiful criterion of the solvability of the diophantine equation $p=x^2+ny^2$. The specific statement is
\thmus{
Let $n$ be an positive integer. Then there is a monic irreducible polynomial $f_n(x)\in\ZZ[x]$ of degree $h(-4n)$ such that if an odd prime $p$ divides neither $n$ nor the discriminant of $f_n(x)$, then $p=x^2+ny^2$ is solvable over $\ZZ$ if and only if $\fracl{-n}{p}{}=1$ and $f_n(x)=0$ is solvable over $\ZZ/p\ZZ$. Here $h(-4n)$ is the class number of primitive positive definite binary forms of discriminant $-4n$. Furthermore, $f_n(x)$ may be taken to be the minimal polynomial of a real algebraic integer $\alpha$ for which $L=K(\alpha)$ is the ring class field of the order $\ZZ[\sqrt{-n}]$ in the imaginary quadratic field $K=\QQ(\sqrt{-n})$.
}
There are some generalizations considering the problem over quadratic fields.
Recently, Harari \cite{bmob} showed that the Brauer-Manin obstruction is the only obstruction for the existence of integral points of a scheme over the ring of integers of a number field, whose generic fiber is a principal homogeneous space of a torus.
After then Dasheng Wei and Fei Xu gave another proof of the result in \cite{multi-norm-tori,multip-type} where the Brouer-Manin obstruction is constructive. This can be used to determine the existence of integral points for the scheme.
Then Dasheng Wei applied the method in \cite{multi-norm-tori} to determine which integers can be written as a sum of two integral squares for some of the quadratic fields $\QQ(\sqrt{\pm p})$ (in \cite{wei1}), $\QQ(\sqrt{-2p})$ (in \cite{wei2}) and so on. That is to determine which integers have integral representations of the form   $x^2+y^2$.

In this paper, we consider  integral representations of the form $ax^2+by^2$ in some other number fields. Although we are interested in certain specific number fields, our formal criterion is for general number fields.

Given a number field $F$ with $\o_F$ its ring of integers. Suppose $a,b$ and $\alpha$ are in $\o_F$ such that $d=-ab$ is not a square in $F$. Let $\XX=\mathbf{Spec}(\o_F[x,y]/(ax^2+by^2-\alpha))$. The general result is:
\thmus{[Theorem \ref{thm_artin_cond}]
Suppose 
one of the following conditions holds:
\enmt{
\it[$(1)$] For every $u$ in $ \o_F^\t$, the equation $x^2-dy^2=u$ is solvable in $\o_F$.
\it[$(2)$] For every  $u\neq1$ in $\o_F^\t$ there exists  a place $\p$ such that the equation $x^2-dy^2=u$ is not solvable in $F_\p$.
}
Then $\XX(\o_F)\neq\emptyset$ if and only if there exists
\eqn{
\prod_{\p\in\Omega_F}(x_\p, y_\p)\in\prod_{\p\in\Omega_F}\XX(\o_{F_\p})
}
such that
\eq{
\psi_{H_L/E}(\tilde f_E(\prod_\p(x_\p,y_\p)))=1.
}}
For notations one can see Section \ref{sec_nota} and for proofs Section \ref{sec_main}.
Due to Wei \cite{wei1,wei2}, the last condition in the theorem is called \emph{Artin condition}.
The theorem shows that if $(1)$ or $(2)$ holds, the Artin condition is the only obstruction to the local-global principle for integral solutions of equation $ax^2+by^2=\alpha$. 

In Section \ref{sec_appl}, we set our focus on some applications of the this theorem.
For a special case ($H_L/F$ is abelian, see Section \ref{sec_abelian}), the Artin condition always holds. Hence the Hasse principal is true.
For the case $F=\QQ$ (Section \ref{sec_q}) and the case $F=\QQ(\sqrt{-p})$ with $a=1,b=q$ 
 where $p$ and $q$ are rational primes  satisfying some conditions
(Section \ref{sec_pq}), condition $(1)$ holds and then we prove that 
the integral local condition in together with the Artin condition completely describe the global integral solvability.
We also give some examples to show the explicit criteria of the solvability.

\section{Solvability by the Artin Condition}\label{sec_artin_cond}

\subsection{Notations}\label{sec_nota}
We fix some notations. Let $F$ be a number field, $\o_F$ the ring of integers of $F$,
$\Omega_F$ the set of all places in $F$ and $\infty_F$ the set of all infinite places in $F$. We write $\p<\infty_F$ for $\p\in\Omega_F\setminus\infty_F$.
Let $F_\p$ be the completion of $F$ at $\p$ and $\o_{F_\p}$ be the valuation ring of $F_\p$ for each $\p<\infty_F$. We also write $\o_{F_\p}=F_\p$ for $\p\in\infty_F$.
The adele ring (resp.  idele group) of $F$ is denoted as $\AA_F$ (resp. $\II_F$).

Let $a,b$ and $\alpha$ be nonzero elements in $\o_F$ and suppose that $d=-ab$ is not a square in $F$. Let $E=F(\sqrt d)$ and $\XX=\mathbf{Spec}(\o_F[x,y]/(ax^2+by^2-\alpha))$ be the affine scheme defined by the equation $ax^2+by^2=\alpha$ over $\o_F$.
The equation $ax^2+by^2=\alpha$ is solvable over $\o_F$ if and only if $\XX(\o_F)\neq\emptyset$.

Denote $R_{E/F}(\GG_m)$ the Weil restriction (see \cite{jmilne-ag}) of $\GG_{m,E}$ to $F$. Let
\eqn{ \varphi: R_{E/F}(\GG_m)\lr\GG_m }
be the homomorphism of algebraic groups which represents
\eqn{ x\mapsto N_{E/F}(x): (E\otimes_FA)^\t\lr A^\t}
for any $F$-algebra $A$. Define the torus $\TT:=\ker\varphi$. 
Let $\XX_F$ be the generic fiber of $\XX$.
Then $\XX_F$ is naturally a $\TT$-torsor by the action of $A$-points:
\aln{
\TT(A)\t \XX_F(A) &\lr \XX_F(A)\\
(u+\sqrt{d}v,x+\frac{\sqrt{d}}{a}y) &\mapsto (u+\sqrt{d}v)(x+\frac{\sqrt{d}}{a}y)
}
It can be seen that 
\eq{\label{eq_in_stab}
\TT(\o_{F_\p})\subseteq \mathbf{Stab}(\XX(\o_{F_\p})):=\set{g\in \TT(F_\p) | g\XX(\o_{F_\p})=\XX(\o_{F_\p})}.
}

Denote by $\lambda$ the embedding of $\TT$ into $R_{E/F}(\GG_m)$.
Tautologically, $\lambda$ induces a natural injective group homomorphism
\eqn{ \lambda_E: \TT(\AA_F)\lr\II_E. }
Let $L=\o_F+\o_F\sqrt{d}$ in $E$ and $L_\p=L\otimes_{\o_F}\o_{F_\p}$  then
\eqn{ \TT(\o_{F_\p})=\set{\beta\in L_\p^\t|N_{E_\p/F_\p}(\beta)=1}. }
It follows that  $\lambda_E(\TT(\o_{F_\p})) \subseteq  L_\p^\t.$
Note that $\lambda_E( \TT(F)) \subseteq E^\t$ in $\II_E$. Hence the following map  induced by $\lambda_E$ is well-defined:
\eqn{
\tilde\lambda_E: \TT(\AA_F)/\TT(F)\TT(\prod_{\p\in \Omega_F} \o_{F_\p})\lr \II_E/E^\t \prod_{\p\in\Omega_F} L_\p^\t.
}
\lemm{\label{lambda_inj}
The map $\tilde\lambda_E$ is injective if one of the following conditions holds:
\enmt{
\it[$(1)$] For every $u$ in $ \o_F^\t$, the equation $x^2-dy^2=u$ is solvable in $\o_F$.
\it[$(2)$] For every  $u\neq1$ in $\o_F^\t$ there exists  a place $\p$ such that the equation $x^2-dy^2=u$ is not solvable in $F_\p$.
}}
\pf{
Recall that $\TT=\ker(R_{E/F}(\GG_m)\lr\GG_m)$.  Therefore we have
\eqn{ \TT(F)=\set{\beta\in E^\t|N_{E/F}(\beta)=1} }
and 
\eqn{ \TT(\o_{F_\p})=\set{\beta\in L_\p^\t|N_{E_\p/F_\p}(\beta)=1}. }
Suppose   $t\in \TT(\AA_F)$ such that $\tilde\lambda_E(t)=1$. Write $t=\beta i$ with $\beta\in E^\t$ and $i\in \prod_\p L_\p^\t$.
Since $t\in \TT(\AA_F)$ we have
\eqn{ N_{E/F}(\beta)N_{E/F}(i)=N_{E/F}(\beta i)=1. }
It follows that
\eqn{ N_{E/F}(\beta)=N_{E/F}(i^{-1})\in F^\t\cap \prod_\p\o_{F_\p}^\t=\o_F^\t. }

If $N_{E/F}(\beta)=N_{E/F}(i)=1$ then we have $\beta\in \TT(F)$ and $i\in\prod_\p\TT(\o_{F_\p})$. Hence $u=\beta i\in \TT(F)\prod_\p\TT(\o_{F_\p})$.

Otherwise $N_{E/F}(\beta)=u$ for some $u\neq1$ in $\o_F^\t$. Then we know that 
 the equation $x^2-dy^2=u$ is  solvable in $F_\p$ for every place $\p$ of $F$. 
This is to say condition $(2)$ does not holds and therefore condition $(1)$ holds. 
It follows that  $x^2-dy^2=u$ is solvable in $\o_F$.
Let  $(x_0,y_0)\in \o_F^2$ be such a solution and let
\eqn{ \zeta=x_0+y_0\sqrt{-1},\gamma=\beta\zeta^{-1}\text{ and }j=i\zeta. }
Then $N_{E/F}(\gamma)=N_{E/F}(j)=1$. Hence $\gamma\in \TT(F)$ and $j\in\prod_\p\TT(\o_{F_\p})$. It follows that $u=\beta i=\gamma j\in \TT(F)\prod_\p\TT(\o_{F_\p})$. 
This finishes the proof.
}

Now we assume that 
\eq{\label{eq_nonempty}
\XX(F)\neq\emptyset.
}
Fixing a rational point $P\in \XX(F)$, since $\XX_F$ is a trivial $\TT$ torsor, we have an isomorphism
\aln{
\phi_P: \prod_{\p\in\Omega_F}\XX(\o_{F_\p}) &\cong \prod_{\p\in\Omega_F}\TT(\o_{F_\p})\\
x &\mapsto P^{-1}x
}
induced by $P$. 
Since we can view $\XX(\o_{F_\p})$ as  an  subset of $\XX(F_\p)$, the composition $f_E:=\lambda_E\phi_P: \prod_\p\XX(\o_{F_\p})\lr  \II_E$ makes sense, mapping 
$x$ to $P^{-1}x$ in $\II_E$. Note that $P$ is in $E^\t\subset \II_E$ since it is a rational point over $F$. It follows that we can 
define the map $\tilde f_E$ to be the composition
\eqn{\xymatrix{
\prod_\p\XX(\o_{F_\p}) \ar@{->}[r]^-{ f_E}  &\II_E       \ar@{->}[r]^-{\t P} &\II_E\\
x                                 \ar@{|->}[r]              &P^{-1}x     \ar@{|->}[r]       &x.
}}
It can be seen that the restriction to $\XX(\o_{F_\p})$ of   $\tilde f_E$ is defined by
\eq{\label{eq_tilde_f_E}
\tilde f_E[(x_\p,y_\p)]=
\cs{
(ax_\p+y_\p\sqrt{d},ax_\p-y_\p\sqrt{d})\in E_{\P_1}\otimes E_{\P_2}	&\text{if }\p\text{ splits in }E/F,\\
ax_\p+y_\p\sqrt{d}  \in E_\P										&\text{otherwise},
}}
where $\P_1$ and $\P_2$ (resp. $\P$) are places of $E$ above $\p$.

Recall that $L=\o_F+\o_F\sqrt{d}$ in $E$ and $L_\p=L\otimes_{\o_F}\o_{F_\p}$. Then $\Xi_L:=E^\t\prod_\p L_\p^\t$ is an open subgroup of $\II_E$.  By the ring class field corresponding to $L$ we mean the class field $H_L$ corresponding to $\Xi_L$ under the class field theory.
Let  $\psi_{H_L/E}: \II_E\lr\Gal(H_L/E)$ be the Artin map. 
For any $\prod_{\p\in\Omega_F}(x_\p, y_\p)\in\prod_{\p\in\Omega_F}\XX(\o_{F_\p})$, noting that $P$ is in $E$, we have 
$\psi_{H_L/E}(f_E(\prod_\p(x_\p,y_\p)))=1$ if and only if
$\psi_{H_L/E}(\tilde f_E(\prod_\p(x_\p,y_\p)))=1$.

\rk{\label{rk_hasse_min}
If  $\prod_{\p\in\Omega_F}\XX(\o_{F_\p})\neq\emptyset$,  then assumption \eqref{eq_nonempty} we made before, that is  $\XX(F)\neq\emptyset$, holds automatically  by Hasse-Minkowski theorem on quadratic equation. Hence We can pick a $F$-point $P$ of $\XX$ and define $\phi_P$.
}

\subsection{The General Theorem}\label{sec_main}
The following is the general theorem we mainly use in this paper, which is motivated by \cite[Corollary 1.4 and 1.6]{multi-norm-tori}.
\thm{\label{thm_artin_cond}
Let symbols be as before. Suppose 
one of the following conditions holds:
\enmt{
\it[$(1)$] For every $u$ in $ \o_F^\t$, the equation $x^2-dy^2=u$ is solvable in $\o_F$.
\it[$(2)$] For every  $u\neq1$ in $\o_F^\t$ there exists  a place $\p$ such that the equation $x^2-dy^2=u$ is not solvable in $F_\p$.
}
Then $\XX(\o_F)\neq\emptyset$ if and only if there exists
\eqn{
\prod_{\p\in\Omega_F}(x_\p, y_\p)\in\prod_{\p\in\Omega_F}\XX(\o_{F_\p})
}
such that
\eq{\label{eq_artin_cond}
\psi_{H_L/E}(\tilde f_E(\prod_\p(x_\p,y_\p)))=1.
}}
\pf{
If $\XX(\o_F)\neq\emptyset$, then 
\eqn{
\tilde f_E\left(\prod_\p\XX(\o_{F_\p})\right)\cap E^\t\prod L_\p^\t\supseteq \tilde f_E(\XX(\o_F))\cap E^\t\neq\emptyset
}
Hence there exists $x\in \prod_{\p\in\Omega_F}\XX(\o_{F_\p})$ such that $\psi_{H_L/E}\tilde f_E(x)=1$.

Conversely, 
we know from Lemma \ref{lambda_inj} that 
\eqn{
\tilde \lambda_E: \TT(\AA_F)/\TT(F)\TT(\prod_\p \o_{F_\p})\lr \II_E/E^{\t}\prod L_\p^\t
}
is injective.
Suppose there exists $x\in \prod_\p\XX(\o_{F_\p})$ such that $\psi_{H_L/E}\tilde f_E(x)=1$ 
(here $\tilde f_E$ makes sense by Remark \ref{rk_hasse_min} in Section \ref{sec_intro}) , i.e. 
$f_E(x)\in E^{\t}\prod L_\p^\t$.
Since $\tilde \lambda_E$ is injective, there are $\tau\in \TT(F)$ and $\sigma\in\TT(\prod_\p \o_{F_\p})$ such that $\tau\sigma=f_E(x)=P^{-1}x$,
i.e. $\tau\sigma(P)=x$. Since $P\in\XX(F)$ and $\TT(\o_{F_\p})\subseteq \mathbf{Stab}(\XX(\o_{F_\p}))$ (see \eqref{eq_in_stab}), it follows that 
\eqn{
\tau(P)=\sigma^{-1}(x)\in \XX(F)\cap \prod_\p\XX(\o_{F_\p})=\XX(\o_F).
}
The proof is complete.
}
The condition \eqref{eq_artin_cond} is called \emph{Artin condition} in Wei's \cite{wei1,wei2}. It interprets the fact that the Brauer-Manin obstruction is the only obstruction for existence of the integral points by conditions in terms of class field theory.
Consequently, the integral local condition in together with the Artin condition completely describe the global integral solvability.
As a result, in the cases where the ring class fields are known it is possible to compute the Artin condition, giving a explicit criteria for the solvability.

\section{Applications}\label{sec_appl}
We now apply the result stated in Theorem   \ref{thm_artin_cond} to proof the following corollaries.

\subsection{On the Cases  $H_L/F$ Abelian}\label{sec_abelian}
Recall that $H_L/F$ is the extension of fields corresponding to $E^\t\prod_\p L_\p^\t$ under the class field theory, where $L=\o_F+\o_F\sqrt{d}$ in $E$ and $L_\p=L\otimes_{\o_F}\o_{F_\p}$.
In the cases that $H_L/F$ is abelian, the Artin condition always holds. Hence the Hasse principal is true.
\thm{\label{thm_trivial_artin_cond}
Suppose $H_L/F$ is abelian, then 
 $\XX(\o_F)\neq\emptyset$ if and only if $\prod_{\p\in\Omega_F}\XX(\o_{F_\p})\neq \emptyset$.
}
\pf{
Since  $H_L/F$ is abelian, the following diagram commutes:
\eqn{\xymatrix{
\II_E \ar@{->}[r]^-{\psi_{H_L/E}} \ar@{->}[d]_{N_{E/F}} &\Gal(H_L/E) \ar@{_(->}[d]\\
\II_F \ar@{->}[r]^-{\psi_{H_L/F}}                       &\Gal(H_L/F)
}}
It follows that
\eqn{
\psi_{H_L/E}(\tilde f_E(\prod_\p(x_\p,y_\p)))
= \psi_{H_L/F}(N_{E/F}( \tilde f_E(\prod_\p(x_\p,y_\p)))) 
= \psi_{H_L/F}(\alpha) = 1
}
where the second equation comes from the definition \eqref{eq_tilde_f_E} of $\tilde f_E$  and
\eqn{
(x_\p+ \sqrt d y_\p)(x_\p- \sqrt d y_\p) =  \alpha\text{ in }E_\P\text{ with }\P\mid \p,
}
while the last equation is obtained by the assumption that $\alpha\in F$.
Hence  the Artin condition always holds and the result follows.
}
\eg{
Let $F=\QQ(\sqrt{-59})$ and $\alpha\in\o_F$. Then the  equation $x^2+2y^2=\alpha$ is solvable over $\o_F$ if and only if
it is solvable over  $\o_{F_\p}$ for all $p\in\Omega_F$.
}
\pf{
In this example, we have $a=1$, $b=2$ and   $d=-2$. Let $E=F(\sqrt{-2})$.
It follows that  $L=\o_F+\o_F\sqrt{-2}=\o_E$ and $H_L=H_E$ the Hilbert field of $E$.
It can be shown that $H_E=EH_F$, and hence $H_L/F$ is abelian. Then the result follows from Theorem \ref{thm_trivial_artin_cond}.
}

\subsection{On the Equation $ax^2+by^2=n$ Over $\QQ$}\label{sec_q}
We consider the case where $F=\QQ$.
\thm{\label{thm_q}
Let $a$, $b$ and $n$ be positive integers. Suppose $d = -ab$ is not a square. Set $E=\QQ(\sqrt d)$,
$L=\ZZ+ \ZZ\sqrt d $ and $H_L$ the ring class field corresponding to $L$. 
Let $\XX=\mathbf{Spec}(\ZZ[x,y]/(ax^2+by^2-n))$.
Then $\XX(\ZZ)\neq\emptyset$ if and only if there exists
\eqn{
\prod_{p\le \infty}(x_p, y_p)\in\prod_{p\le \infty}\XX(\ZZ_p)
}
such that
\eqn{
\psi_{H_L/E}(\tilde f_E(\prod_p(x_p,y_p)))=1
}
where  $\tilde f_E$ is defined the same as in \eqref{eq_tilde_f_E} 
except $F=\QQ$.
}
\pf{
Since $d=-ab<0$ it is clear that $x^2-dy^2=-1$ is not solvable over $\RR$, which is to 
say the condition $(2)$ of Theorem \ref{thm_artin_cond} holds since the only units of $\ZZ$ are $\{\pm1\}$. 
Then the Proposition applies, whence the result follows.
}
We now give an example where explicit  criterion is obtained using this result.
\eg{
Let $n$ be a positive integer and $f(x)=x^4-x^3+x+1\in\ZZ[x]$. 
Write $n=2^{s_1}\t7^{s_2}\t\prod_{k=1}^g p_k^{e_k}$ and define
$D=\{p_1,p_2,\dots,p_g\}$ and 
\aln{
D_1&=\set{p\in D| \fracn{-14}{p}=1\text{ and }f(x)\mod p\text{ irreducible}},\\
D_2&=\set{p\in D| \fracn{-14}{p}=1\text{ and }f(x)\text{ splits into two irreducible factors}}.
}
Then the diophantine equation $2x^2+7y^2=n$ is solvable over $\ZZ$ if and only if
\enmt{
\itm{1} $n\t 2^{-s_1}\equiv\pm1\pmod8$,
\itm{2} $\fracn{n\t 7^{-s_2}}{7}=1$,
\itm{3} for all $p\nmid2\t7$, $\fracn{-14}{p}=1$ for odd $e, v_p(n)=e$,
\itm{4} and
\eqn{\cs{
\sum{p_i\in D_1}e_i\equiv0\pmod2             &\text{ if }D_1\neq\emptyset,\\
1+s_1+s_2+\sum_{p_j\in D_2}e_j\equiv0\pmod2  &\text{ if }D_1=\emptyset.
}}}}
\pf{
In this example, we have $a=2$, $b=7$,$d=-14$ and $d\equiv2\pmod4$. Let $E=\QQ(\sqrt{-14})$.
It follows that  $L=\ZZ+\ZZ\sqrt{-14}=\o_E$ and $H_L=H_F=E(\alpha)$ the Hilbert field of $E$ 
where the minimal polynomial of $\alpha$ is $f(x)$. The Galois group $\Gal(H_L/E)=<\sqrt{-1}>\cong\ZZ/4\ZZ$.
Let $\XX=\mathbf{Spec}(\ZZ[x,y]/(2x^2+7y^2-n))$ and
\eqn{
\tilde f_E[(x_p,y_p)]=
\cs{
(2x_p+y_p\sqrt{-14},ax_p-y_p\sqrt{-14})	&\text{if }p\text{ splits in }E/\QQ,\\
2x_p+y_p\sqrt{-14}                      &\text{otherwise}.
}}
Then by Theorem \ref {thm_q}, $\XX(\ZZ)\neq\emptyset$ if and only if there exists
\eqn{
\prod_{p\le \infty}(x_p, y_p)\in\prod_{p\le \infty}\XX(\ZZ_p)
}
such that
\eqn{
\psi_{H_L/E}(\tilde f_E(\prod_p(x_p,y_p)))=1.
}
Next we compute these conditions in details.
By  a simple  computation the local condition  \eqn{
\prod_{p\le \infty}\XX(\ZZ_p)\neq\emptyset
} is equivalent to 
\eq{\label{eq_eg_q_local}
\cs{
n\t 2^{-s_1}\equiv\pm1\pmod8,\\
\fracn{n\t 7^{-s_2}}{7}=1,\\
\text{for all }p\nmid2\t7, \fracn{-14}{p}=1\text{ for odd }e, v_p(n)=e.
}}
For Artin condition, let $(x_p,y_p)_p\in\prod_p\XX(\ZZ_p)$. Then 
\eq{\label{eq_local_decomp}
(x_p+ \sqrt {-14} y_p)(x_p- \sqrt {-14}  y_p) =  n\text{ in }E_\P\text{ with }\P\mid p.
}
And since $H_L/E$ is unramified, for any $p$ we have 
\eq{\label{eq_psi_1}
1=\cs{
   \psi_{H_L/E}(p_\P)\psi_{H_L/E}(p_{\bar\P}),    &\text{ if }p=\P\bar\P\text{ splits in }E/\QQ,\\
   \psi_{H_L/E}(p_\P),                             &\text{ if }p=\P\text{ inert in }E/\QQ,
}}
where $p_\P$ (resp. $p_{\bar\P}$)  is in $\II_E$ such that its $\P$ (resp. $\bar\P$) component is $p$ and the other
components are $1$.
We calculate $\tilde f_E[(x_p,y_p)]$ separately:
\enmt{
\it If $p=2$, $2=\P_2^2$ in $E/\QQ$. Suppose $\P_2=\pi_2\o_{E_{\P_2}}$ for $\pi_2\in\o_{E_{\P_2}}$.
Since $\P_2^2$ is principal in $E$ but $\P_2$ is not, we have $\psi_{H_L/E}((\pi_2)_{\P_2})=-1$.
By \eqref{eq_local_decomp} we have
\eqn{ v_{\P_2}(2x_2+\sqrt{-14}y_2)=v_{\P_2}(2x_2-\sqrt{-14}y_2)=\frac{1}{2}v_{\P_2}(2n)=v_2(2n)=s_1+1.}
It follows that \aln{ 
\psi_{H_L/E}(\tilde f_E[(x_2,y_2)])&=\psi_{H_L/E}((2x_2+\sqrt{-14}y_2)_{\p_2})\\
                                   &=(-1)^{v_{\P_2}(2x_2+\sqrt{-14}y_2)}=(-1)^{s_1+1},
}where $\tilde f_E[(x_2,y_2)]$ is also regarded as an element in $\II_E$ such that the components above $2$ are given by the value 
of $\tilde f_E[(x_2,y_2)]$ and $1$ otherwise.
\it If $p=7$, a similar argument shows that $\psi_{H_L/E}(\tilde f_e[(x_2,y_2)])=(-1)^{s_2}$.
\it If $\fracn{-14}{p}=1$  then by \eqref{eq_psi_1} we can distinguish the following cases:
\enmt{[(i)]
    \it $f(x)\mod p$ splits into linear factors.
    Then $\psi_{H_L/E}(p_\P)=\psi_{H_L/E}(p_{\bar\P})=1$ and  $\psi_{H_L/E}(\tilde f_E[(x_p,y_p)])=1$.
    \it $f(x)\mod p$ splits into two irreducible factors. 
    Then $\psi_{H_L/E}(p_\P)=\psi_{H_L/E}(p_{\bar\P})=-1$.
    It follows that \aln{ 
    \psi_{H_L/E}(\tilde f_E[(x_p,y_p)])&=\psi_{H_L/E}((2x_p+\sqrt{-14}y_p)_\P)(2x_p-\sqrt{-14}y_p)_{\bar\P})\\
                                       &=(-1)^{v_\P(2x_p+\sqrt{-14}y_p)+v_{\bar\P}(2x_p-\sqrt{-14}y_p)}=(-1)^e,
    }where $v_p(n)=e$ since \aln{
    v_\P(2x_p+\sqrt{-14}y_p)&+v_{\bar\P}(2x_p-\sqrt{-14}y_p)\\
                            &=v_p(2x_p+\sqrt{-14}y_p)+v_p(2x_p-\sqrt{-14}y_p)=v_p(n).
    }
    \it $f(x)\mod p$  irreducible.
    Then $\psi_{H_L/E}(p_\P)=-\psi_{H_L/E}(p_{\bar\P})=\pm\sqrt{-1}$.
    It follows that \aln{ 
    \psi_{H_L/E}&(\tilde f_E[(x_p,y_p)])=\psi_{H_L/E}((2x_p+\sqrt{-14}y_p)_\P)(2x_p-\sqrt{-14}y_p)_{\bar\P})\\
                &=(\pm\sqrt{-1})^{v_\P(2x_p+\sqrt{-14}y_p)+v_{\bar\P}(2x_p-\sqrt{-14}y_p)}
                (-1)^{v_{\bar\P}(2x_p-\sqrt{-14}y_p)}\\
                &=(\pm\sqrt{-1})^e(-1)^a
    }where $v_p(n)=e$  and $a=v_p(2x_p-\sqrt{-14}y_p)$ (in $\QQ_p$, $0\le a\le e$). By Hensel lemma, we can choose the local solution
    $(x_p,y_p)$ suitably, such that $a$ can rich any value between $0$ and $e$. 
    Hence $ \psi_{H_L/E}(\tilde f_E[(x_p,y_p)]) =\pm(\sqrt{-1})^e$ with the sign chosen freely.
}
\it If $\fracn{-14}{p}=-1$ then $p$ inert in $E/\QQ$. By \eqref{eq_psi_1} we have $\psi_{H_L/E}(\tilde f_E[(x_p,y_p)])=1$.
\it At last if $p=\infty$, since $H_L/E$ is unramified, we have $\psi_{H_L/E}(\tilde f_E[(x_\infty,y_\infty)])=1$.
}
Putting the above argument together, we know the Artin condition is 
\eq{\label{eq_eg_q_artin}
\cs{
\sum{p_i\in D_1}e_i\equiv0\pmod2             &\text{ if }D_1\neq\emptyset,\\
1+s_1+s_2+\sum_{p_j\in D_2}e_j\equiv0\pmod2  &\text{ if }D_1=\emptyset.
}}
The proof is done if we put the local condition \eqref{eq_eg_q_local} and the Artin condition \eqref{eq_eg_q_artin} together.
}

\subsection{On the Equation $x^2+qy^2=\alpha$ Over $\QQ(\sqrt{-p})$}\label{sec_pq}
At last we consider the  representations of the form $x^2+qy^2$ over a class of imaginary quadratic fields $\QQ(\sqrt{-p})$ where
\eq{\label{eq_pq_cond}
p\text{ and }q\text{ are rational primes of the form }4k+1\text{ satisfying }\fracn{p}{q}=-1\text{ or }\fracl{p}{q}{4}=-1.
}
\thm{\label{thm_pq}
Let $p$, $q$ satisfy \eqref{eq_pq_cond}  and $F=\QQ(\sqrt{-p})$ be quadratic field.
Let $\alpha$ be an element in $\o_F$. Set $d = -q$, $E=F(\sqrt d)$,
$L=\o_F+ \o_F\sqrt d $ and $H_L$ the ring class field corresponding to $L$. 
Let $\XX=\mathbf{Spec}(\o_F[x,y]/(x^2+qy^2-\alpha))$.
Then $\XX(\o_F)\neq\emptyset$ if and only if there exists
\eqn{
\prod_{\p\in\Omega_F}(x_\p, y_\p)\in\prod_{\p\in\Omega_F}\XX(\o_{F_\p})
}
such that
\eqn{
\psi_{H_L/E}(\tilde f_E(\prod_\p(x_\p,y_\p)))=1
}
where  $\tilde f_E$ is defined the same as in \eqref{eq_tilde_f_E} 
except $a=1$ and $b=q$.
}
\pf{
Suppose that  \eqref{eq_pq_cond} holds. Since $p$ and $q$ are distinct and of the form $4k+1$, we have $\o_F=\ZZ+\ZZ\sqrt{-p}$ and $q$ is not a square in $F$. We claim that the condition $(1)$ of Theorem \ref{thm_artin_cond} holds.
Since $\o_F^\t={\pm 1}$, it's enough to prove that $x^2+qy^2=-1$ has an integral solution over $\o_F$. 
Since $p$ and $q$ are distinct and of the form $4k+1$, we have $\o_F=\ZZ+\ZZ\sqrt{-p}$ and $q$ is not a square in $F$.
If $\fracn{p}{q}=-1$, or $\fracn{p}{q}=1$ and $\fracl{p}{q}{4}=\fracl{q}{p}{4}=-1$, then $x^2-pqy^2=-1$ is solvable over $\ZZ$ (\cite[p. 228]{dirichlet}); otherwise $\fracn{p}{q}=1$, $\fracl{p}{q}{4}=-1$ and $\fracl{q}{p}{4}=1$, so $x^2-pqy^2=p$ is solvable over $\ZZ$ (\cite[Corollary 1.4]{wei_diophantine}).
Both cases indicate that there exists $x_0,y_0\in\o_F$ such that $x_0^2+qy_0^2=-1$.
Hence Theorem \ref{thm_artin_cond} yields the result.
}

\newnoindbf{{Acknowledgment}
The authors would like to thank Dasheng Wei  and Congjun Liu for many helpful discussions and comments.
}

\bibliography{\thefilename}
\bibliographystyle{amsplain} 
\end{document}